\documentclass[10pt,reqno]{amsart}
\usepackage[T1]{fontenc}
\usepackage[latin1]{inputenc}
\usepackage{amsmath}
\usepackage{amssymb}
\usepackage{amsthm}
\usepackage{amsbsy}
\usepackage{enumerate}
\usepackage{graphics}
\usepackage[normalem]{ulem}
\usepackage{tikz}
\usetikzlibrary{patterns}
\usepackage{xcolor}
\usepackage[most]{tcolorbox}
\tcbuselibrary{breakable, skins}
\usepackage{caption}
\usepackage{subfigure}
\usepackage[a4paper,
left=2cm, right=2cm,
top=3cm, bottom=3cm]{geometry}

\usepackage{hyperref}
\hypersetup{
    colorlinks,
    linkcolor={luh-dark-blue},
    citecolor={luh-dark-blue},
    urlcolor={orange}
}

\definecolor{luh-dark-blue}{rgb}{0.0, 0.313, 0.608}
\definecolor{navyblue}{rgb}{0.0,0.0,0.5}
\definecolor{zaffre}{rgb}{0.0, 0.08, 0.66}
\definecolor{white}{rgb}{1.0, 1.0, 1.0}
\definecolor{darkblue}{rgb}{0.0, 0.2, 0.6}
\definecolor{darkgray}{rgb}{0.66, 0.66, 0.66}
\definecolor{lightgray}{rgb}{0.83, 0.83, 0.83}

\newtheorem{satz}{Proposition}[section]
\newtheorem{lem}[satz]{Lemma} 
\newtheorem{remark}[satz]{Remark}
\newtheorem{thm}[satz]{Theorem}

\newtheorem{cor}[satz]{Corollary}

\newcommand{\chookrightarrow}{\mathrel{\lhook\joi34nrel\relbar\kern-.8ex\joinrel\lhook\joinrel\rightarrow}}

\newcommand{\R}{\mathbb{R}}

\newcommand{\N}{\mathbb{N}}		   
\newcommand{\Z}{\mathbb{Z}}

\newcommand{\T}{\mathbb{T}}

\newcommand{\e}{\varepsilon}

\normalsize
\DeclareMathOperator{\supp}{supp}
\DeclareMathOperator{\esssup}{ess-sup}

\numberwithin{equation}{section}

	\title[Gardner-Ostrovsky type equation]{Existence, regularity, and symmetry of periodic traveling waves for Gardner--Ostrovsky type equations}

\author{Gabriele Bruell}
\address{Centre for Mathematical Sciences, Lund University, 22100 Lund, Sweden}
\email{gabriele.brull@math.lth.se}
\author{Long Pei}
\address{School of Mathematics (Zhuhai), Sun Yat-sen University, 519082 Zhuhai, China}
\email{peilong@mail.sysu.edu.cn}
\thanks{Date: \today }

\begin{document}

	\maketitle
	\allowdisplaybreaks

\begin{abstract}
We study the existence, regularity, and symmetry of periodic traveling solutions to a class of Gardner--Ostrovsky type equations, including the classical Gardner--Ostrovsky equation, the (modified) Ostrovsky, and the reduced (modified) Ostrovsky equation. The modified Ostrovsky equation is also known as the short pulse equation. The Garder--Ostrovsky equation is a model for internal ocean waves of large amplitude. We prove the existence of nontrivial, periodic traveling wave solutions using local bifurcation theory, where the wave speed serves as the bifurcation parameter.  Moreover,  we give a regularity analysis for periodic traveling solutions in the presence as well as absence of the Boussinesq dispersion. We see that the presence of Boussinesq dispersion implies smoothness of periodic traveling wave solutions, while its absence may lead to singularities in the form of peaks or cusps. Eventually, we study the symmetry of periodic traveling solutions by the method of moving planes. A novel feature of the symmetry results in the absence of Boussinesq dispersion is that we do not need to impose a traditional monotonicity condition or a recently developed reflection criterion on the wave profiles to prove the statement on the symmetry of periodic traveling waves.
 
\end{abstract}

\vspace{10pt}


\vspace{10pt}

\section{Introduction}
We consider the Gardner--Ostrovsky type equation of the form
\begin{equation}\label{eq:GardnerOstrovsky}
	(u_{t}+n(u)_x+\beta u_{xxx})_{x}=\gamma u, \quad \quad  \beta,\gamma \geq 0.
\end{equation}
 The function $u=u(t,x)$ is the unknown and depends on  $(t,x)\in \R_+\times \R$. The nonlinearity $n(\cdot)$ is a nontrivial, real analytic function with $n(0)=n^{\prime}(0)=0$. The subscripts denote the partial derivatives of $u$ with respect to the subscript.  For 
\begin{equation}\label{quadratic-cubic nonlinearity}
n(u)=\frac{\sigma}{2}u^2+\frac{\alpha}{3}u^3,\qquad \sigma, \alpha \in \R,
\end{equation}
equation \eqref{eq:GardnerOstrovsky} is  referred to as the Gardner--Ostrovsky equation or the extended rotation-modified Korteweg--de Vries equation, which describes  long internal waves of large
amplitude \cite{HPT}. The parameter $\beta$ represents the Boussinesq (small scale) dispersion due to the non-hydrostaticity caused by the finite fluid depth, while the parameter $\gamma$ represents the Coriolis (large scale) dispersion due to the Earth's rotation \cite{MR3771099}.  The cubic nonlinearity appears in internal wave dynamics when the traditional quadratic, hydrodynamic nonlinearity is anomalously small \cite{HPT}, or when approximately large-amplitude waves are considered \cite{MiBa,MR2184887}.  For certain choices of the parameters $\alpha, \beta$ and $ \gamma$, equation \eqref{eq:GardnerOstrovsky} reduces to different well-known equations. For instance, if $\gamma=0$ and $ \sigma=1$, \eqref{eq:GardnerOstrovsky} with nonlinearity \eqref{quadratic-cubic nonlinearity} becomes the Gardner equation
\[
u_{t}+uu_{x}+\alpha u^{2}u_{x}+\beta u_{xxx}=0,
\]
which then reduces to the famous Korteweg--de Vries (KdV) equation provided $\alpha=0$. If  $\gamma >0$, 
$ \sigma=1$ and $\alpha=0$, then \eqref{eq:GardnerOstrovsky} becomes the Ostrovsky equation (also known as the rotation-modified KdV equation)
\begin{equation}\label{eq:O}
	(u_{t}+uu_{x}+\beta u_{xxx})_{x}=\gamma u
\end{equation}
and the reduced Ostrovsky equation if also $\beta=0$. Instead, when $\gamma>0$, $\alpha=1$ and $\sigma=0$, then \eqref{eq:GardnerOstrovsky} becomes the modified Ostrovsky equation (also known as short pulse equation)
\begin{equation}\label{eq:mO}
(u_{t}+u^2u_{x}+\beta u_{xxx})_{x}=\gamma u,
\end{equation}
and the reduced modified Ostrovsky equation (reduced short pulse equation), if also $\beta=0$. 
\medskip

\subsection*{Related results concerning the Ostrovsky and modified Ostrovsky equation}
In the case when $\beta \neq 0$,
the Ostrovsky equation is locally well-posed in $H^s(\R)$ for $s\geq -\frac{3}{4}$ \cite{Tsugawa,MR3805879} and ill-posed for $s<-\frac{3}{4}$ \cite{Tsugawa}.
Local well-posedness for the modified Ostrovsky equation is confirmed  in Sobolev spaces 
with regularity index $s\geq \frac{1}{4}$, while ill-posedness holds for $s<\frac{1}{4}$ (see \cite{MR3740577,MR3799685}). 
{\color{black} Solitary waves are constructed via variational arguments as ground state solutions in $H^{1}(\R)$-based energy spaces  (see \cite{MR4299849,MR2291873}) or  $H^{2}(\R)$-based energy spaces (see \cite{MR4094633}), or by Fenichel singular perturbation theory \cite{MR2578800} with exponential decay.  The  stability is also considered in related function spaces in the these results (the stability of solitary solutions   in $H^{1}(\R)$-based  spaces is also studied in \cite{MR2995236} for generalized Ostrovsky equations with homogeneous nonlinearities).  In contrast, the non-existence of  solitary solutions is confirmed for  $\beta\gamma>0$ (see \cite{Galkin, Leonov}). It is worth to mention that the  multi-pulse traveling waves are obtained via  geometric singular perturbation theory and Lyapunov-Schmidt reduction for the modified Ostrovsky equation in \cite{MR2575365}.}  
In \cite{MR3614681}, periodic traveling waves in $H^{2}(\T)$  are constructed by viewing the (modified) Ostrovsky equation as a perturbation of the (modified) KdV via the Lyapunov-Schmidt reduction argument. In addition, the linear stability of those periodic solutions is confirmed, and spectral stability can be obtained when the Coriolis force $\gamma$ is relatively small \cite{MR3614681}. 
In \cite{MR3527630}, the orbital stability of periodic traveling wave solutions for the (modified) Ostrovsky equation is proved.

\medskip

In the case when $\beta=0$, that is, when the Boussinesq dispersion is absent, equations \eqref{eq:O} and \eqref{eq:mO} are often referred to as the reduced (modified) Ostrovsky equation. 
Wave-breaking for solutions with smooth initial data  for the reduced Ostrovsky equation is confirmed in \cite{MR2684307}.
Smooth periodic solutions exist for $c\in (0,c^{*})$ where $c^{*}$ is uniquely determined in \cite{MR3685174} (see also \cite{BD} for a generalization). In addition, when $c=c^{*}$ the solutions are peaked, and an explicit form of the peaked solution in terms of polynomials is given for $c^{*}=\frac{\pi^2}{9}$. The peaked periodic waves are shown to be spectrally unstable in the space of square-integrable periodic functions with zero mean and the same period in \cite{MR4163826}. 
The linear instability and uniqueness of the peaked periodic wave are confirmed for the reduced Ostrovsky equation in \cite{MR3936896}.  Furthermore, in \cite{ MR3094592,MR2993126} the authors show that the peaked periodic traveling wave solution for the reduced (modified) Ostrovsky equation is in some sense between wave breaking and global well-posedness. More precisely, it is shown in  \cite{ MR3094592,MR2993126} that solutions corresponding to smooth initial data $u_0$ satisfying a certain sign condition exist globally in time, while the solution exhibits wave breaking when the sign condition is not definite but sign-changing. It turns out that the peaked solutions satisfy the sign condition almost everywhere except at the peaks, where the condition is not well-defined.

\medskip

Local well-posedness for the Gardner--Ostrovsky type equation \eqref{eq:GardnerOstrovsky} with nonlinearity \eqref{quadratic-cubic nonlinearity} is studied in \cite{MR4083153}. To the best of our knowledge, there are so far no existence results in the literature concerning periodic traveling waves  for \eqref{quadratic-cubic nonlinearity} with non-zero $\alpha$ and $\sigma$. 

\medskip

\subsection*{Main results}
The present paper aims to study the existence, regularity and symmetry of periodic traveling wave solutions to \eqref{eq:GardnerOstrovsky}. In our studies, we {\color{black} assume that $\beta\geq 0$ and fix $\gamma=1$  for convenience\footnote{{\color{black} Note that the results still hold for other $\gamma>0$ after proper adjustment in the formulation.}}}. Moreover, $n$ is assumed to be real analytic with $n(0)=n^\prime(0)=0$. 
We prove the existence of even, periodic traveling wave solutions for \eqref{eq:GardnerOstrovsky} utilizing local bifurcation theory. Furthermore,  we conduct an {a priori} regularity analysis on periodic traveling solutions for the Gardner--Ostrovsky type equation.  We show that the presence of  Boussinesq dispersion $\beta>0$ implies that periodic traveling waves are smooth. On the other hand, when the Boussinesq dispersion is absent, that is $\beta=0$, we prove that small amplitude solutions are still smooth, while large amplitude solutions may exhibit singularities in the form of peaks or cusps. This result can be seen in the frame of the often observed dependence of the regularity of a traveling solution on its amplitude and propagation speed {\color{black}(see \cite{MR4002168,MR4504564,MR4257625,MR4458407} for uni-directional models and \cite{MR3902471} for bidirectional model)}. Eventually, we investigate the symmetry of periodic traveling wave solutions and apply our results to the Gardner--Ostrovsky equation, which is \eqref{eq:GardnerOstrovsky} with nonlinearity \eqref{quadratic-cubic nonlinearity}.
A summary of our results reads:

\medskip

\begin{tcolorbox}

\emph{{Existence}}\\
There exist continuous, even, $2\pi$-periodic traveling wave solutions for \eqref{eq:GardnerOstrovsky} (Theorem \ref{thm:local}).

\medskip

\emph{{Regularity}}\\
Let $u(t,x)=\phi(x-ct)$ be a continuous, $2\pi$-periodic solution of  \eqref{eq:GardnerOstrovsky} and set $F(\phi)=-c\phi + n(\phi)$.
\begin{itemize}
	\item If $\beta>0$, then $\phi$ is smooth (Theorem \ref{prop:beta}).
	\item If $\beta =0$ and $F^\prime(\phi)<0$, then $\phi$ is smooth, i.e., small amplitude traveling wave solutions are smooth. If $F^\prime(\phi)\leq 0$ and $F^\prime(\phi)(\bar x)=0$ for some $x\in \T$, then $\phi$ exhibits a singularity in the form of a peak or a cusp 	(Theorem \ref{prop:beta0}). More precisely, for $|x-\bar x|\ll 1$ we have that
	\[
	|\phi(x)-\phi(\bar x)|\eqsim |x-\bar x|^\frac{2}{a}, \qquad {\color{black}a:=\min \{n\geq 2 \mid F^{(n)}(\bar \phi)\neq 0\}}.
	\]

\end{itemize}

\medskip

\emph{{Symmetry}}\\
Let $u(t,x)=\phi(x-ct)$ be a continuous, $2\pi$-periodic  solution of \eqref{eq:GardnerOstrovsky}.
\begin{itemize}
	\item If $\beta>0$ and $\phi$ has a single crest per period, then $\phi$ is symmetric (Theorem \ref{thm:symmetry_beta_0}).
	\item If $\beta =0$ and $\phi$ has a unique global maximum and minimum per period, then $\phi$ is symmetric and has a single crest per period
	(Theorem \ref{thm:symmetry_beta_1}).
\end{itemize}

\end{tcolorbox}

\medspace

The motivation for focusing on periodic traveling waves comes from the fact that solitary traveling waves cannot appear for $\beta\gamma>0$ (see Theorem \ref{Anti-soliton theorem}).  This fact is referred to as the \emph{anti-soliton} theorem for the Ostrovsky equation \cite{Galkin, Leonov}.

\medskip

	The study of symmetries of solutions in partial differential equations has a long history. A classical argument for {\color{black} this issue} is the method of moving planes, which goes back to Alexandrov \cite{Alk} and was refined by the works of Serrin \cite{Serr}, Gidas, Ni, and Nirenberg \cite{GNN} and many others in various contexts. The strong connection between symmetry and traveling waves has also been observed in the context of water waves, where the first results are by Garabedian \cite{Garabedian} for periodic traveling waves and by Craig and Sternberg \cite{Craig} for solitary solutions. The symmetry property arises from a structural condition related to the dispersion relation of the equation, see also \cite{MR3603270, 2021Symmetry,MR4113186}.
{\color{black}The symmetry  of} periodic traveling waves is more involved compared with the case of solitary solutions due to the lack of decay at infinity. To overcome this issue, classical symmetry results for periodic traveling waves are based on a monotonicity condition, which is often imposed between a unique trough and crest period. 
Recently, a \emph{reflection criterion} which does not require any local monotonicity assumption, is put forward in \cite{2021Symmetry} for periodic traveling solutions to a class of nonlocal dispersive equations with very weak dispersion. {\color{black}We focus on periodic traveling waves in this manuscript.}
When the Boussinesq dispersion is absent ($\beta=0$) in \eqref{eq:GardnerOstrovsky},  we  {\color{black} confirm the symmetry}  if the wave is merely assumed to have a unique global maximum and minimum per period. Thus, for $\beta=0$, we can confirm the symmetry of solutions without imposing either a monotonicity condition or the recently established reflection criterion on the wave profile. On the other hand, when the Boussinesq dispersion is involved, we prove the symmetry of traveling solutions with a single crest per period (i.e., if the so-called global monotonicity condition is satisfied) 

\medskip

\medskip

\medskip

We conclude this section with an outline of the paper. In Sect. \ref{Preliminaries}, we collect some preliminaries including an anti-soliton theorem (Theorem \ref{Anti-soliton theorem}) for the Gardner--Ostrovsky type equation  \eqref{eq:GardnerOstrovsky} when $\gamma\beta>0$, which excludes the existence of solitary solutions and motivates our focus on periodic traveling waves.  
	 In Sect. \ref{sect:existence and regularity}
	 we prove the existence result for even, periodic traveling wave solutions by means of local bifurcation theory. Sect. \ref{sect:regularity analysis} is devoted to a regularity analysis of periodic traveling  wave solutions for \eqref{eq:GardnerOstrovsky} with and without  Boussinesq dispersion.
In Sect. \ref{symmetry of smooth solutions}, we prove the symmetry results of both smooth traveling waves and waves with peaks/cusps. Eventually, in Sect. \ref{sect:application to GO}, we apply the above results to the Gardner--Ostrovsky equation, which is \eqref{eq:GardnerOstrovsky} with nonlinearity \eqref{quadratic-cubic nonlinearity}.

\bigskip

\section{Preliminaries}\label{Preliminaries}
Let us consider the Gardner--Ostrovsky type equation 
\begin{equation}\label{eq:GardnerOstrovsky1}
	(u_{t}+n(u)_x+\beta u_{xxx})_{x}=\gamma u, \quad \quad \beta, \gamma \in \R.
\end{equation}

It is \emph{a priori} clear from the equation that for $\gamma \neq 0$ any sufficiently smooth, spatially periodic or integrable solution of \eqref{eq:GardnerOstrovsky1} has zero mean, that is
\[
\int u (t,x)\,dx =0,
\]
where the integral is respectively taken over one period or the real line.

\subsection{Dispersion relation and \emph{anti-soliton theorem}}

Solving the linearized {\color{black}Gardner}--Ostrovsky type equation reveals a relation between the wave frequency $\omega$ and the wave number $k$ given by
\begin{equation}\label{eq:dispersion}
	\omega(k) =\frac{\gamma}{k}-\beta k^3.
\end{equation}
This relation is called the \emph{dispersion relation}. The phase velocity $c_p$ and the group velocity $c_g$ can be computed from the dispersion relation as follows:
\begin{equation}\label{eq:c}
	c_p(k):= \frac{\omega(k)}{k}= \frac{\gamma}{k^2}-\beta k^2\qquad \mbox{and}\qquad c_g(k):=\frac{d}{dk}\omega(k)=-\frac{\gamma}{k^2}-3\beta k^2.
\end{equation}
In Figure \ref{F:d}, we plot the phase velocity $c_p$ with $\gamma>0$ for the linearized Gardner--Ostrovsky equation for both  $\beta>0$ (solid line) and $\beta <0$ (dashed line).

\begin{center}
	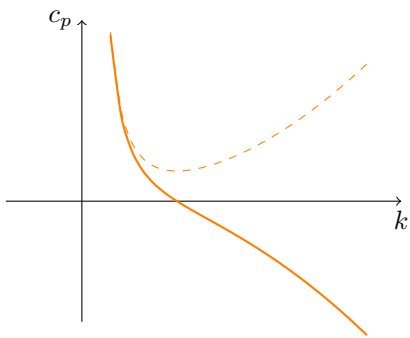
\begin{figure}[h!]
\begin{tikzpicture}[yscale=0.8]
	\draw[->] (-1, 0) -- (4.2, 0) node[below] {$k$};
	\draw[->] (0, -2) -- (0, 3) node[left] {$c_p$};
	\draw[scale=0.5, xscale=2.5, domain=0.3:3, smooth, variable=\x, orange, thick] plot ({\x}, {0.5/(\x*\x)-0.5*(\x*\x)});	
	\draw[scale=0.5, xscale=2.5, domain=0.3:3, smooth, variable=\x, orange, dashed] plot ({\x}, {0.5/(\x*\x)+0.5*(\x*\x)});
\end{tikzpicture}
\caption{Plot of the phase velocity $c_p$ for $k>0$ for the Gardner--Ostrovsky equation in the case $\gamma>0$ and $\beta>0$ (solid line) as well as $\beta <0$ (dashed line)} \label{F:d}
\end{figure}
\end{center}

The present study focuses on the existence, regularity and symmetry of traveling waves for the case $\gamma>0$ and $\beta\geq 0$. 
We start our observation by an \emph{anit-soliton theorem}, stating that no solitary traveling wave solutions exist for the Gardner--Ostrovsky type equation when $\beta \gamma>0$.  If $u(t,x)=\phi(x-ct)$ is a traveling wave solution of \eqref{eq:GardnerOstrovsky1} with wave speed $c>0$, then $\phi$ solves the ordinary fourth-order differential equation
\begin{equation}\label{eq:ODE}
(-c\phi + n(\phi)+\beta \phi_{xx})_{xx}=\gamma \phi.
\end{equation}
We say that $\phi$ is a smooth solitary solution of \eqref{eq:ODE} if $\phi$ solves \eqref{eq:ODE} and $|\phi(x)|\to 0$ for $|x|\to \infty$ together with its derivatives. 

\begin{thm}[Anti-soliton theorem]\label{Anti-soliton theorem}
If $\beta \gamma>0$, no nontrivial smooth solitary solutions exist for \eqref{eq:ODE}.
	\end{thm}

\begin{proof}
Observe that \eqref{eq:ODE} can be written as
\begin{equation}\label{eq:ode2}
\phi^{\prime\prime}= \frac{1}{\beta} \left( c\phi -n(\phi)+\gamma v\right),\qquad v^{\prime\prime}=\phi
\end{equation}
and the {\color{black}Hamiltonian}
\begin{equation}\label{eq:energy E}
E(\phi,v):= \frac{1}{2} \left( (\phi^\prime)^2+\frac{\gamma}{\beta}(v^\prime)^2\right)  + \frac{1}{\beta} \left(-\frac{c}{2}\phi^2+N(\phi)-\gamma v\phi \right),
\end{equation}
where $N^\prime(s):=n(s)$,
is constant along solutions of \eqref{eq:ode2}. If $(\phi,c)$ solves \eqref{eq:ode2}, then 
\begin{align*}
\frac{d}{dx}E(\phi,v) &= \phi^\prime \phi^{\prime\prime} + \frac{\gamma}{\beta} v^\prime v^{\prime\prime} + \frac{1}{\beta} \left(-c\phi \phi^\prime + n(\phi)\phi^\prime - \gamma v\phi^\prime - \gamma v^\prime \phi\right)= \left(\phi^{\prime\prime} -\frac{1}{\beta}\left( c\phi -n(\phi)+\gamma v\right)\right) \phi^\prime =0.
\end{align*}

Therefore, $E(\phi,v)=0$ if $\phi$ is a nontrivial solitary solution of \eqref{eq:ODE} and $v^{\prime\prime}=\phi$. 
Note that any smooth solution of \eqref{eq:ODE}  fulfills a priori a zero mean condition. Therefore, if a global maximum of $\phi$ is attained at $\bar x\in \R$, then there exists on open neighborhood $\Omega$ of $\bar x$, with $\phi>0$ in $\Omega$ and $\phi=0$ on $\partial \Omega$. Given \eqref{eq:ode2}, we find that $v$ is convex on $\Omega$. Note that   $\phi(x_0)=0$ for $x_0\in \partial \Omega$ and $\frac{\gamma}{\beta}>0$. Then, in view of the fact  $E(\phi, v)$ is constantly zero, we can use  \eqref{eq:energy E} to get
\[
\phi^\prime(x_0)=v^\prime(x_0)=0 \qquad \mbox{for} \quad x_0\in \partial \Omega.
\]
But this is a contradiction to the Hopf boundary lemma, which states that a convex function $v$ on $\Omega$ takes its maximal value at the boundary with a strictly positive outward derivative unless $v$ is constant. This implies by $v^{\prime\prime}=\phi$ that $\phi=0$ on $\Omega$. We conclude that the global maximum of $\phi$ is zero, which is a contradiction to $\phi$ being nontrivial with zero means.

\end{proof}

\begin{remark} \emph{We would like to emphasize that for the Ostrovsky equation, which has the form of \eqref{eq:GardnerOstrovsky} with $\alpha=0$ in \eqref{quadratic-cubic nonlinearity}, the anti-soliton theorem was already known and proved for $\gamma \beta >0$ \cite{Galkin, Leonov}. {\color{black} The proof here and in  \cite{Galkin, Leonov} all make use of the Hamiltonian of the corresponding equation in consideration. The difference is that the argument in \cite{Galkin} is based on the linearized equation, while the argument in \cite{Leonov} is based on expansion of the solution   and the analysis of coefficients. 
} Furthermore,  the authors in \cite{Zhang} investigate symmetry and uniqueness of solitary solutions when $\gamma \beta <0$. 
	}
\end{remark}

The above anti-soliton theorem motivates our focus on periodic traveling waves. 

\subsection{Functional analytic setting}
We denote by $\T:=\R / 2\pi \Z$  the one-dimensional torus, which is identified with $[0,2\pi) \subset \R$. 
Let  $\mathcal D(\T)=C^\infty(\T)$ and denote by $\mathcal{S}(\Z)$ the space of rapidly decaying functions. Then the (periodic) Fourier transform $\mathcal{F}:\mathcal{D}(\T)\to \mathcal S(\Z)$ is defined by
\[
(\mathcal{F}f)(k)=\hat f(k):=\frac{1}{{2\pi}}\int_{\T}f(x)e^{-ixk}\,dx,
\]
and any $f\in \mathcal{D}(\T)$ can be written as
\[
f(x)= \sum_{k\in \Z}\hat f(k)e^{ixk}.
\]
By duality, the Fourier transform extends uniquely to $\mathcal{F}:\mathcal{D}^\prime(\T)\to \mathcal S^\prime(\Z)$.
Here $\mathcal D^\prime(\T)$ and $\mathcal{S}^\prime(\T)$ are the dual spaces of $\mathcal{D}(\T)$ and $\mathcal{S}(\T)$, respectively.
We say a function $f:\T \to \R$ belongs to the space $L^p(\T)$, $1\leq p<\infty$, if and only if
\[
\|f\|_{L_p}^p:=\int_\T |f|^p(x)\, dx < \infty,
\]
and $f\in L^\infty(\T)$ if and only if $\|f\|_\infty:= \esssup_{x\in \T}|f(x)|<\infty$. 

We now introduce the periodic Zygmund spaces on which we perform our subsequent analysis.
Let $(\varphi)_{j\geq 0}\subset C_c^\infty(\R)$ be a family of smooth, compactly supported functions satisfying
\[
\supp \varphi_0 \subset [-2,2],\qquad \supp \varphi_j \subset [-2^{j+1},-2^{j-1}]\cup [2^{j-1},2^{j+1}] \quad\mbox{ for}\quad j\geq 1,
\]
\[
\sum_{j\geq 0}\varphi_j(\xi)=1\qquad\mbox{for all}\quad \xi\in\R,
\]
and for any $n\in\N$, there exists a constant $c_n>0$ such that 
\[\sup_{j\geq 0}2^{jn}\|\varphi^{(n)}_j\|_\infty\leq c_n.\]

For $s>0$, the periodic Zygmund space denoted by $ \mathcal{C}^s(\T)$ consists of functions $f$ satisfying
\[
\|f\|_{\mathcal{C}^s(\T)}:=\sup_{j\geq 0}2^{sj}\left\|\sum_{k\in \Z} e^{ik(\cdot)}  \varphi_j(k)\hat f(k)\right\|_\infty < \infty.
\]

Eventually, for $\alpha \in (0,1)$, we denote by $C^\alpha(\T)$ the space of $\alpha$-H\"older continuous functions on $\T$. 
If $k\in \N$ and $\alpha\in (0,1)$, then $C^{k,\alpha}(\T)$ denotes the space of $k$-times continuously differentiable functions whose $k$-th derivative is $\alpha$-H\"older continuous on $\T$ {\color{black}with the norm
\begin{equation*}
\|u\|_{C^{k,\alpha}}:=\sum_{l\leq k}\sup_{x\in \T}|D^{l}u(x)|+\sup_{x,y\in \T,  x\neq y}\left\{\frac{|u(x)-u(y)|}{|x-y|^{\alpha}}\right\}.
\end{equation*}
}
To lighten the notation we write $C^s(\T)=C^{\left \lfloor{s}\right \rfloor, s- \left \lfloor{s}\right \rfloor }(\T)$ for $s\geq 0$. 
As a consequence of Littlewood--Paley theory, we have the relation $\mathcal{C}^s(\T)=C^s(\T)$ for any $s>0$ with $s\notin \N$; that is, the H\"older spaces on the torus are completely characterized by Fourier series. If $s\in \N$, then $C^s(\T)$ is a proper subset of $\mathcal{C}^s(\T)$ and
\[
C^1(\T)\subsetneq C^{1-}(\T)\subsetneq \mathcal{C}^1(\T).
\]
Here, $C^{1-}(\T)$ denotes the space of Lipschitz continuous functions on $\T$. For more details, we refer to \cite[Chapter 13]{T3}. We denote by $X_{0}$ the subspace of functions of mean zero in a function space $X$. In particular, we use in the following  $\mathcal{C}_{0}^{r}(\T)$ and $C^s_{0}(\T)$  to denote the subspace of periodic functions with zero mean {\color{black}in the Zygmund space $\mathcal{C}^{r}(\T)$ and  H{\"o}lder space $C^s(\T)$ for $s\geq 0$, respectively.} Moreover, we add the subscript \emph{even} whenever we restrict further to the subspace of even functions.

\subsection{Fourier multipliers}

A  Fourier multiplier on $\Z$ is a complex-valued function that defines a linear operator $L$ via multiplication on the Fourier side, that is
\[
\widehat{Lf}(k)=m(k)\hat f(k).
\] 
The function $m$ is also called the symbol of the multiplier operator $L$. 

If $f$ belongs to $C^{a}(\T)$ for $a>\frac{1}{2}$ and is real-valued, then $f$ has a Fourier series expression of the form
\[
f(x)=\sum_{k\in \Z} f_k e^{ikx},\qquad \mbox{where}\quad f_k=f_{-k},
\]
 and the Fourier series is convergent. The Fourier coefficients $(f_k)_{k\in \Z}$ are given by
\[
f_k = \hat f(k)=\int_\T f(x)e^{-ixk}\,dx.
\]
Note that the subset of functions $f$ in $C^{a}(\T)$ for $a>\frac{1}{2}$ with $\hat f(k)=\hat f(-k)$ for all $k\in Z$ is invariant under the nonlinearity $n$.

\bigskip

\section{Existence of periodic traveling waves}\label{sect:existence and regularity}

We consider the Gardner--Ostrovsky type equation of the form
\begin{equation}\label{eq:Gardner-Ostrovsky}
(u_{t}+n(u)_x+\beta u_{xxx})_{x}=u, \quad  \quad \beta\geq 0,
\end{equation}
where $n$ is nontrivial, real analytic with $n(0)=n^\prime(0)=0$.
Here, we have fixed $\gamma=1$ and $\beta \geq 0$. 
With the ansatz $u(t,x)=\phi(x-ct)$, we get the following 
 steady equation
 \begin{equation}\label{eq:eqs}
 	(-c\phi + n(\phi)+\beta \phi_{xx})_{xx}=\phi.
 \end{equation}
 
 Setting $D=i\partial_x$ and introducing the nonlocal Fourier multiplier operators $D^{-2}$ and $L_{\beta, c}:=\frac{1}{c +\beta D^{2}}$, it will be convenient to rewrite \eqref{eq:eqs} in the nonlocal form
 \begin{equation}\label{equation}
 	-\phi + D^{-2}L_{\beta, c}\phi + L_{\beta, c} n(\phi)=B_\phi,
 \end{equation}
where $B_\phi\in \R$ is an integration constant. Since the operator $D^{-2}$ is a Fourier multiplier operator which is only well-defined on functions of zero mean, the integration constant is determined by
 \[
 B_\phi = \frac{1}{c}\int_\T n(\phi)\,dx.
 \]
 
 \medspace
 
 \emph{
In what follows, we say that \emph{$u(t,x)=\phi(x-ct)$ is a periodic traveling wave solution} of the Gardner--Ostrovsky type equation \eqref{eq:Gardner-Ostrovsky}, if $\phi$ is a continuous, $2\pi$-periodic solution of \eqref{eq:eqs} with zero mean.
}

\medspace
\begin{lem}[Properties of the linear operators]\label{regularity of operators}
	Set $s_\beta = 2$ for any $\beta >0$ and $s_0=0$.
	The operator $L_{\beta, c}$ is a (smoothing) operator of order $s_\beta$ and the operator $D^{-2}$ is a smoothing operator of order 2, that is
	\[
	L_{\beta, c}:\mathcal{C}_0^{r}(\T)\to \mathcal{C}_0^{r+s_\beta}(\T)\qquad \mbox{and}\qquad 	D^{-2}:\mathcal{C}_0^{r}(\T)\to \mathcal{C}_0^{r+2}(\T)
	\]
	are bounded operators for any $r\geq 0$.  Furthermore, the action of  the operator $D^{-2}$ can be expressed as a convolution with the kernel function $K$ as
	\[
 D^{-2}\phi =K*\phi,\qquad \mbox{where}\quad K(x)=\frac{1}{4\pi}(|x|-\pi)^2-\frac{\pi}{12},\qquad x \in [-\pi,\pi].
	\]
If $\beta >0$, the action of the operator $L_{\beta, c}$ can be expressed as a convolution with the kernel function $G_{\beta, c}$ as
	\[
		L_{\beta, c}\phi = G_{\beta, c}*\phi, \qquad \mbox{where}\qquad G_{\beta,c}(x)=\frac{1}{2\pi c}+\frac{\beta \pi}{4\pi c\sinh(\pi \frac{c}{\beta})}\cosh((\pi-|x|)\frac{c}{\beta})-\frac{\beta^2}{4\pi c^2}.
	\]
	\end{lem}

\begin{proof}
	The smoothing properties of $L_{\beta,c}$ and $D^{-2}$ are shown in e.g. \cite[Proposition 2.78]{BCD} and \cite[Proposition 3.1]{BD}, respectively, while the explicit expression for the convolution kernels $K$ and $G_{\beta, c}$ follow from
\[
K(x)=\frac{1}{\pi} \sum_{k=1}^\infty |k|^{-2}\cos(xk)=\frac{1}{4\pi}(|x|-\pi)^2-\frac{\pi}{12}
\]
and
\[
G_{\beta, c}(x)=\frac{1}{\pi} \sum_{k=1}^\infty \frac{1}{c+\beta k^2}\cos(xk)=\frac{1}{2\pi c}+\frac{\beta \pi}{4\pi c\sinh(\pi \frac{c}{\beta})}\cosh((\pi-|x|)\frac{c}{\beta})-\frac{\beta^2}{4\pi c^2}
\]
for $x\in[-\pi,\pi]$, where we used that the symbols of the nonlocal operators are the Fourier coefficients of the kernel functions, respectively. The explicit expressions of the series can be found in \cite{PBM}.
\end{proof}

\begin{remark}\label{remark}\emph{
	If $\phi$ is continuous with zero mean on $\T$, then $D^{-2}\phi$ belongs to the class $C^2(\T)$. Note that $K^\prime(x)=\frac{1}{2\pi}x -\frac{ x}{2|x|}$. Then, we have
	\begin{align*}
	[K*\phi]^\prime(x)&= \frac{1}{2\pi}\int_\T (x-y)\phi(y)\,dy -\frac{1}{2}\int_\T \frac{(x-y)}{|x-y|}\phi(y)\,dy= -\frac{1}{2\pi}\int_\T y\phi(y)\,dy -\frac{1}{2} \int_{-\pi}^x\phi(y)\,dy+\frac{1}{2}\int_{x}^\pi \phi(y)\,dy.
	\end{align*}
 Taking the second derivative of $K*\phi$  yields that
\[
	[K*\phi]^{\prime\prime}(x)=- \phi(x).
\]
Thus, $D^{-2}\phi$  belongs to $C^2(\T)$ with $\partial_x^2 D^{-2}\phi=-\phi$.
}
\end{remark}

We now apply analytic bifurcation theory to confirm the existence of small amplitude solutions for \eqref{equation}. 
 We are going to {\color{black}work on $C^{\alpha}_{\textrm{even}}(\T)$, $\alpha\in (1,2)$, the restriction to even functions in $C^{\alpha}(\T)$.} Define
\begin{equation}\label{Definition of nonlinearity F}
G(c,\phi):=-\phi + D^{-2}L_{\beta, c}\phi + L_{\beta, c} n(\phi)-B_{\phi}.
\end{equation}
Note that $G:\R_+\times C^a_{0,\textrm{even}}(\T)\to C^a_{0,\textrm{even}}(\T)$ is {\color{black}analytic} and $G(c,0)=0$ for all $c\in \R_+$. Thus, the set $\{(c,0)\mid c\in \R_+\}\subset \R_+\times C^a_{0,\textrm{even}}(\T)$ constitutes a trivial solution line.
The Frech\'et derivative of $G$ with respect to $\phi$ is given by
\begin{equation*}
D_{\phi}G(c,\phi)[\psi]=-\psi + \frac{1}{c D^{2}+\beta D^{4}}\psi + \frac{1}{c +\beta D^{2}}n^\prime(\phi)\psi
\end{equation*}
so that 
\begin{equation}\label{eq: linearization at phi equal 0}
D_{\phi}G(c,0)[\psi]=-\psi + \frac{1}{c D^{2}+\beta D^{4}}\psi,
\end{equation}
since $n^\prime(0)=0$.
\begin{lem}\label{lemma:Fredholm operator}
$D_{\phi}G(c,0)$ is a Fredholm operator of index zero. Moreover, if $\beta \in [0,1)$, then there exists $c_k:=\frac{1}{k^2}-\beta k^2>0$ such that the kernel of $D_{\phi}G(c_k,0)$ is one-dimensional and  spanned by $\{\cos(k\cdot)\}$.
\end{lem}
\begin{proof}
Recall that the operator $D^{-2}L_{\beta, c}$ {is bounded linear operator}  form $C^{a}_{0,\textrm{even}}(\T)$ to $C^{a+2+s_\beta}_{0,\textrm{even}}(\T)$,
the latter  being compactly embedded into $C^{a}_{0,\textrm{even}}(\T)$. Therefore,  $D^{-2}L_{\beta, c}=\frac{1}{cD^{2}+\beta D^4}$ is a compact operator on $C^a_{0,\textrm{even}}(\T)$. It is then clear that $D_{\phi}G(c,0)$ is a compact perturbation of the identity operator. Thereby, it is a Fredholm operator of index zero.

In order to show that the kernel of $D_{\phi}G(c,0)$ is one-dimensional, we use the Fourier transform and observe that $\psi \in \ker D_{\phi}G(c,0)$ if and only if
\[
\left(1-\tfrac{1}{ck^2+\beta k^4}\right)\hat \psi (k) =0\qquad \mbox{for all}\qquad k\in \Z\setminus\{0\}.
\]
Then $\hat \psi$ has nontrivial support if and only if
\begin{equation}\label{c}
	c=\frac{1}{k^{2}}-\beta k^{2}, \quad k\in \Z\setminus\{0\}.
\end{equation}
For $\beta\geq 0$, the right-hand side is strictly decreasing. Therefore,  the dimension of the $\ker D_\phi G(c,0)$ is at most one, and exactly one if \eqref{c} is satisfied. Since $\beta \in [0,1)$, we find at least one value $c>0$ satisfying  \eqref{c}.  The kernel of $D_{\phi}G(c,0)$ is then spanned by $\{\cos(k\cdot)\}$ because $\hat \psi(k)=-\hat \psi(-k)$.

\end{proof}

The existence of a local curve of nontrivial solutions bifurcating from the zero solution then follows from Lemma \ref{lemma:Fredholm operator} and the {\color{black}Crandall-Rabinowitz} theorem \cite{MR1956130,Kielhofer-book} if the transversality condition 
\[
D^2_{c\phi}G(c_k,0)[\psi] \notin \mbox{{\color{black}ran}}\, D_\phi G(c_k,0)\qquad \mbox{for}\quad \psi \in \ker D_\phi G(c_k,0)
\]
is true.
Recall that
\[
D_\phi G(c,\phi)=-\mbox{Id} + D^{-2}L_{c,\beta}+ L_{c,\beta}n^\prime(\phi),
\]
where $L_{c,\beta}= (c+\beta D^2)^{-1}$. The derivative of $L_{c,\beta}$ with respect to $c$ is given by
\[
\frac{d}{dc}L_{c,\beta} = \frac{d}{dc}(c+\beta D^2)^{-1} =-(c+\beta D^2)^{-1}\left(\frac{d}{dc}(c+\beta D^2)\right)(c+\beta D^2)^{-1}=-L_{c,\beta}^2.
\]
Hence, we have 
 \[
D^2_{c,\phi}G(c_k,0)=-D^{-2}L_{c,\beta}^2 =\frac{D^{-2}}{(c+\beta D^2)^2}.
\]
Assume for a contradiction that there exists $f\in C^{a}_{0,\textrm{even}}(\T)$ such that
\[
D^2_{c,\phi}G(c_k,0)[\psi]= D_\phi G(c_k,0)f \qquad \mbox{for}\quad \psi \in \ker D_\phi G(c_k,0).
\]
Since $ \ker D_\phi G(c_k,0)=\{\cos(k\cdot)\}$, where $k$ is such that $c_k=\frac{1}{k^2}-\beta k^2$, we find that $\hat f(j)=0$ for all $j\in \N\setminus\{k\}$ and
\[
-\frac{2}{k^2(c_k+\beta k^2)^2}=\left(1-\frac{1}{c_k k^2+\beta k^4}\right)\hat f(k),
\]
which is a contradiction since the right-hand side vanishes due to the definition of $c_k$, while the left-hand side is strictly negative. We conclude that the transversality condition is fulfilled. With the above preparation, the following theorem on the existence of periodic traveling solutions is a direct consequence of the Crandall--Rabinowitz theorem .

\begin{thm}[Local bifurcation] \label{thm:local}
	Let $a\in (1,2)$ and $\beta \in [0,1)$. There exists a local  bifurcation curve of nontrivial $2\pi$-periodic real-valued solutions of \eqref{equation} emanating from the bifurcation point $(c_*,0)$, where $c_*=1-\beta$. The curve is given by
	\[
	\{(c(\e),\phi(\e))\mid |\e|<\e_0\}\subset \R_+ \times C^a_{0,\mathrm{even}}(\T).
	\]
	{\color{black}In addition, $c$ and $\phi$ are analytic functions on $(-\e,\e)$.}
\end{thm}

{\color{black}
\begin{remark}
	\emph{
We remark that the period is chosen to be $2\pi$ for convenience. The bifurcation can be done for waves with other period and the parameter $\beta$ is connected to the wavelength.   It is also worth to mention that the authors in \cite{MR4176870} consider the bifurcation diagram of the capillary-gravity Whitham equation. The steady form of this equation is rewritten also in the way that a nonlocal operator is coupled with the nonlinearity. The continuation of local bifurcation is given there, and a larger solution space by two-parameter bifurcation is performed when the symbol of the operator is not monotone on the positive half real line. The latter may  be connected to the case here when $\beta<0$.
}x
\end{remark}
}

\bigskip

\section{Regularity of periodic traveling waves}\label{sect:regularity analysis}

In this section, we state our results for the steady equation
\begin{equation}\label{eq:rGOt}
	F(\phi)+\beta \phi_{xx}=-D^{-2}\phi + B_\phi,
\end{equation}
where the function $F$ satisfies the following assumptions.

\textbf{Assumption ($F$)}
\begin{itemize}
	\item[(i)] $F:\R \to \R$ is real analytic;
	\item[(ii)] $F(0)=0$; 
	\item[(iii)]  $F^\prime(0)=-c$ for some $c>0$ and $S:= \{ \phi \in \R \mid F^\prime(\phi)=0\}\neq \emptyset$.
\end{itemize}
The set $S$ contains all critical points of $F$. In what follows, we denote (if they exist) by $\bar \phi_-<0$ and $\bar \phi_+>0$ the two critical values which are closest to the origin, respectively. Notice that at least one of them exists by assumption (F)(iii). When it does not confuse, we  use the simple notation 
\begin{equation}\label{eq:bar}
	\bar \phi \in \{\bar \phi_-,\bar \phi_+\},
\end{equation} whenever $\bar \phi_-$ {\color{black}or} $\bar \phi_+$ exists.

\begin{remark}\emph{
The steady Gardner--Ostrovsky type equation is covered by the  more general form of  \eqref{eq:rGOt} as
\[
F(\phi)=-c\phi + n(\phi).
\]	
In particular, the assumption (F) is satisfied in view of $n$ being nontrivial, real analytic with $n(0)=n^\prime(0)=0$.
}
\end{remark}

\subsection{The case of non-zero Boussinesq dispersion $\boldsymbol{\beta >0}$}

We observe that the presence of Boussinesq dispersion $\beta>0$ implies that any continuous solution of \eqref{eq:rGOt} is \emph{a priori} smooth.

\begin{satz}[Regularity of solutions for $\beta >0$]\label{prop:beta}
	Let  $\phi$ be a continuous $2\pi$-periodic solution of \eqref{eq:rGOt} for $\beta>0$. Then $\phi$ is smooth.
\end{satz}
\begin{proof}
	By assumption (F), equation \eqref{eq:rGOt} can be written as
\begin{equation}\label{equation regularity analysis 1}
\phi= D^{-2}L_{\beta, c}\phi + L_{\beta, c} r(\phi)+B_\phi,
\end{equation}
where $r(\phi)=\sum_{n=2}^\infty \frac{F^{(n)}(0)}{n!}\phi^n$. 
Given Lemma \ref{regularity of operators},  the operator
\begin{equation}
\phi\mapsto D^{-2}L_{\beta, c}\phi + L_{\beta, c} r(\phi)
\end{equation}
maps a continuous function $\phi$ into the space $\mathcal{C}_0^{r+2}(\T)$. Then \eqref{equation regularity analysis 1} implies that $\phi\in\mathcal{C}_0^{r+2}(\T)$. Repeating the above procedure, we conclude that $\phi\in C^{\infty}(\T)$.
\end{proof}

\subsection{The case of zero Boussinesq dispersion $\boldsymbol{\beta =0}$ } 

If $\beta=0$, equation \eqref{eq:rGOt} reduces to the steady equation for the \emph{reduced Gardner--Ostrovsky} type equation, that is
\begin{equation}\label{eq:RGO}
	F(\phi)+D^{-2}\phi=B_\phi
\end{equation}
 In this case, traveling wave solutions are not \emph{a priori} smooth but might exhibit singularities in the form of peaks or cusps if their amplitude is large enough. 

\bigskip

\begin{satz}[Regularity of solutions for $\beta =0$] \label{prop:beta0}
	Let $\phi$ be a continuous $2\pi$-periodic solution of \eqref{eq:rGOt}. Then, the following hold:
	\begin{itemize}
		\item[(a)]$\phi$ is smooth on any open set where $F^\prime(\phi)<0$;
		\item[(b)] If  $\phi$ has an isolated local maximum or minimum at some $\bar x$  with $\phi(\bar x)=\bar \phi$, where $\bar \phi$ is as in \eqref{eq:bar},
		then $\phi$ does not belong to  $C^1(\T)$. More precisely, let $a:=\min \{n\geq 2 \mid F^{(n)}(\bar \phi)\neq 0\}$, then $\phi$ is precisely $\frac{2}{a}$-H\"older continuous at $\bar x$ and
		\begin{equation}\label{eq:H}
		|\phi(x)-\bar \phi| \eqsim \left(\frac{a!}{2F^{(a)}(\bar \phi)}\bar \phi \right)^\frac{1}{a}|x-\bar x|^{\frac{2}{a}}\qquad \mbox{for}\quad |x-\bar x|\ll 1.
		\end{equation}
	
	\end{itemize}
\end{satz}

\begin{proof}
	\begin{itemize}
		\item[(a)]   {\color{black}The proof is analog to the proof of Lemma 3.3 in \cite{MR4504564} which in turn is based on \cite{MR4002168,MR3902471}, and the analysis of $F$ follows from Lemma 4.1 in \cite{MR3902471}.} We include details here for convenience and completeness.
		Let $\phi$ be a continuous, $2\pi$-periodic solution of \eqref{eq:rGOt}.  
		If  $F^\prime(\phi)<0$ on some open set $\Omega$, then for any point $x\in \Omega$, there exists a  neighborhood $B_{\e}(x)$ with $\e>0$ sufficiently small such that $F^\prime(\phi)<0
		$ on $B_{\e}(x)$.  Restricting such $\phi$ on $B_{\e}(x)$, we find that  the map $F:\phi\mapsto F(\phi)$ is strictly  decreasing from the set  $\{\phi \in C(\T)\,| \,F^\prime(\phi)<0\}$ to the set of $2\pi$-periodic continuous functions. Hence, by inverse function theorem, there exists a smooth inverse of $F$ denoted by $F^{-1}$ such that 
		\begin{equation}\label{inverse F on the equation}
			\phi(y)=F^{-1} (F(\phi))(y)=F^{-1}(-D^{-2}\phi+B_\phi)(y), \quad \quad y\in B_{\e}(x).
		\end{equation}
		Note that for any $2\pi$-periodic continuous function $\phi$, we have $-D^{-2}\phi\in  C_{0}^{2}(\T)$ (cf. Remark \ref{remark}). Therefore, the right-hand side of \eqref{inverse F on the equation} is twice continuously differentiable so that $\phi$  is also continuously differentiable on $B_{\e}(x)$. Repeating this procedure, we find that $\phi$ is smooth on $\Omega$. 
		\item[(b)] Assume that $\phi$ has a local maximum or minimum at $\bar x$  with $\bar \phi:=\phi(\bar x)$, where $\bar \phi$ as in \eqref{eq:bar} so that  $F^\prime(\bar \phi)=0$. In view of (F)(i), we can write \eqref{eq:rGOt} as
		\[
		\sum_{n\geq 2}^\infty \frac{F^{(n)}(\bar \phi)}{n!}(\phi-\bar \phi)^n = -D^{-2}\phi + B_\phi -F(\bar \phi).
		\]
		By continuity of $\phi$ we have that the right-hand side $G:=-D^{-2}\phi + B_\phi -F(\bar \phi)$ above is twice continuously differentiable and $G(\bar x)=G^\prime(\bar x)=0$, since $\phi$ has a local maximum or minimum at $x=\bar x$. Applying Taylor expansion around $\bar x$ on $G$, we find that
		\[
			\sum_{n\geq 2}^\infty \frac{F^{(n)}(\bar \phi)}{n!}(\phi-\bar \phi)^n(x) = \frac{1}{2}\phi (\xi)(x-\bar x)^2,\qquad \mbox{for some}\quad \xi \in [\bar x, x].
		\]
Using that $a=\min \{n\geq 2 \mid F^{(n)}(\bar \phi)\neq 0\}$ one obtains that
		\begin{equation}\label{eq:regu}
	\lim_{x\to \bar x}  \frac{|\phi(x)-\bar \phi|^a}{|x-\bar x|^2}= \frac{a!}{2|F^{(a)}(\bar \phi)|}|\bar \phi|.
		\end{equation}
		Suppose that  $\phi$ were continuously differentiable at $\bar{x}$.  Then
		\[
		0= 	\lim_{x\to \bar x}   \left(\frac{|\phi(x)-\bar \phi|}{|x-\bar x|}\right)^a(x-\bar x)^{a-2}= \frac{a!}{2|F^{(a)}(\bar \phi)|}|\bar \phi|,
		\]
		since $\phi$ takes a local extremum at $\bar x$. But this is a contradiction because $|\bar \phi|\neq 0$. From \eqref{eq:regu} we obtain that $\phi$ is precisely $\frac{2}{a}$-H\"older continuous at $\bar x$, and   \eqref{eq:H} follows  from \eqref{eq:regu}.
	\end{itemize}
\end{proof}

\begin{remark}[Bound on the amplitude if $a$ is even]\label{rem:amplitude} \emph{If  $a=\min \{n\geq 2 \mid F^{(n)}(\bar \phi)\neq 0\}$ is even, then any $2\pi$-periodic solution of \eqref{eq:rGOt} is bounded from above by $\bar \phi_+$ and from below by $\bar \phi_-$, whenever $\bar \phi_+$ and/ or $\bar \phi_-$ exist. This follows immediately from
		\begin{equation}\label{eq:rem}
	\sum_{n\geq a}^\infty \frac{F^{(n)}(\bar \phi)}{n!}(\phi-\bar \phi)^n = -D^{-2}\phi + B_\phi -F(\bar \phi)
	\end{equation}
	and the fact that the right-hand side of \eqref{eq:rem} is strictly convex (concave) on any connected subset $P^+$ of $\T$ where $\phi$ is positive (negative). Suppose $\bar \phi_+>0$ exists and $\phi$ is a solution of \eqref{eq:rGOt} with  $\phi(\bar x)=\bar \phi_+$ and $\max \phi=\phi(x_M)>\bar \phi_+$. Without loss of generality, we may assume that $\bar x$ and $x_M$ are inner points of a connected subset $P^+\subset\T$ where $\phi$ is positive. Since the right-hand side of \eqref{eq:rem} is strictly convex on $P^+$, we deduce that it has at most one local extremum (which needs to be a minimum) on $P^+$. However, since $a$ is even, the left-hand side has a local minimum at $x=\bar x$ and a local maximum at $x=x_M$, which yields a contradiction. {\color{black}Note that it could be quite difficult to get a similar concise result when $a$ is odd . We will study  the bound of solutions for an odd $a$, the convexity of the  extreme waves, and the nodal properties in a future work.}
}
\end{remark}

\bigskip

\section{Symmetry of traveling waves}\label{symmetry of smooth solutions}

As in the previous section, we consider the equation
\begin{equation}\label{eq:GOs}
F(\phi)+\beta \phi_{xx}=-D^{-2}\phi +B_\phi,
\end{equation}
where $F$ satisfies assumption (F).

\medspace

By Lemma \ref{regularity of operators}, we know that the action of the operator $D^{-2}$ on functions with zero mean can be written as a convolution 
\[
D^{-2}\phi = K*\phi,
\]
where $K(x)=2\sum_{k=1}^\infty |k|^{-2}\cos(xk)$.
Note that $K$ is even, continuous, smooth outside the origin, strictly decreasing on $(0,\pi)$ and belongs to the Sobolev space $ W^{1,\infty}(\T)$. {\color{black}The touching lemma and boundary point lemmas below will be used in the proof of symmetry of periodic traveling waves. They follow  previous arguments in   \cite{2021Symmetry}, but modification   is needed based on the property of $F$. So, we include the proofs for completeness. }

\begin{lem}[Touching lemma within one period]\label{touching lemma}
Let $\phi, \tilde \phi$ be two $2\pi $-periodic continuous solutions of \eqref{eq:GOs} with $B_\phi = B_{\tilde \phi}$. Suppose that  $\phi \geq \tilde{\phi}$ on $[\lambda,\lambda +\pi ]$ and $\phi-\tilde{\phi}$ is odd with respect to $\lambda$.  Then, 
\begin{itemize}
\item either $\phi=\tilde{\phi}$  on $\T$, or 
\item $\phi>\tilde{\phi}$ on $(\lambda,\lambda+\pi)$.
\end{itemize}
\end{lem}
\begin{proof}
{\color{black}
Set $w:=\phi-\tilde \phi$. Then $w$ is  $2\pi $-periodic with $w(x)\geq 0$ on $(\lambda, \lambda+\pi )$, and   $w(2\lambda-x)=-w(x)$.   The function $w$ solves the equation
\begin{equation}\label{eq:sup-sub solution equation1}
	F(\phi)-F(\tilde \phi)+\beta w_{xx}+K*w= 0.
	\end{equation}

Let us assume that $w$ is not identically zero, but there exists a point $\bar x \in (\lambda, \lambda+\pi )$ such that $w(\bar x)=0$. Recall that  $w$ is smooth (cf. Proposition \ref{prop:beta}) if $\beta \neq 0$. In addition, $\beta w_{xx}(\bar x)\geq 0$ because  $w$ attains a local minimum at $\bar x$. Then, evaluating  \eqref{eq:sup-sub solution equation1}   at $\bar x$ yields that
\begin{equation}\label{L2 w at xbar}
K*w(\bar x) \leq 0.
\end{equation}
Denoting by $G(y):=K(\bar x-y)-K(\bar x+y-2\lambda) $, we split the integral above in the following way:
\begin{align}\label{eq:L1 on w at xbar}
\begin{split}
K*w(\bar x)=&\int_{\lambda}^{\bar x}G(y)w(y)\,dy+\int_{\bar x}^{\lambda +\pi }G(y)w(y)\,dy,
\end{split}
\end{align}
where $\lambda<\bar x <\lambda+\pi $.
We aim to show that $G(y)$ is positive on $(\lambda, \lambda+\pi)$, and this will indicate that $K*w(\bar x)$ is positive via \eqref{eq:L1 on w at xbar}
 if $w$ is not identically zero. Then a  contradiction with \eqref{L2 w at xbar} appears.
Note that  $G$ is $2\pi $-periodic and odd with respect to $\lambda$ and
$
 G(\lambda)=G(\lambda+\pi )=0
$. 
%
Set $z=y-\bar x $ and $v=2(\bar x-\lambda)$. In view of the monotonicity of  $K$ on $(0,\pi )$, it is clear that  $G(y)=0$  if and only if $v\in 2\pi \Z$ or $v \in -2z+2\pi \Z$. On one hand, we have $v=2(\bar x-\lambda)\in (0,2\pi )$ so that $v\notin 2\pi \Z$.  On the other hand, we observe that $v \in -2z+2\pi \Z$ if and only if there exists $n\in \Z$ such that
$y=\lambda+2\pi n \notin (\lambda,\lambda+\pi)$. Therefore, we have  $G(y)\neq 0$ on $(\lambda,\lambda+\pi )$. The nonnegativity of $G(y)$ on $(\lambda,\lambda+\pi )$ then implies $G(y)> 0$ on this interval.  
Recalling \eqref{eq:L1 on w at xbar} and $w(y)\geq 0$ on $(\lambda, \lambda+\pi )$, we obtain that 
$K*w(\bar x)>0$
unless $w$ is identically zero. Hence, we get  a contradiction to \eqref{L2 w at xbar}, and conclude that $\phi>\tilde \phi$ on $(\lambda,\lambda+\pi)$ unless $\phi=\tilde \phi$ on $\T$.
 }
\end{proof}

While the touching lemma is related to a strong maximum principle, the following lemma plays a role as the Hopf boundary point lemma does for elliptic equations.

\begin{lem}[Boundary point lemma]\label{boundary lemma} Let $\phi, \tilde \phi$ be two $2\pi $-periodic continuous solutions of \eqref{eq:GOs} with $B_\phi = B_{\tilde \phi}$. Suppose that $\phi \geq \tilde \phi$ on $[\lambda, \lambda + \pi ]$, $\phi-\tilde \phi$ is odd with respect to $\lambda$, and $\phi-\tilde \phi$ is continuously differentiable in a neighborhood of $\lambda$. Then either
\begin{itemize}
\item $\phi = \tilde \phi$ on $\T$, or
\item $ (\phi-\tilde \phi) ^\prime(\lambda)> 0$.
\end{itemize}

\end{lem}

\begin{proof} Set $w(x):=\phi(x)-\bar{\phi}(x)$. Since $w$ is odd with respect to $\lambda$, we have that $w(\lambda)=0$ and $w^\prime(\lambda)\geq 0$ in view of $\phi \geq \tilde \phi$ on $[\lambda,\lambda+\pi]$. Now, let us assume that $w$ is not identically zero and assume for a contradiction that $w^\prime(\lambda)=0$. Using assumption (F), we have that $w$ solves 
\begin{equation}\label{eq:equation for boundary point lemma 0}
	-cw+r(\phi)-r(\bar\phi)+\beta w_{xx}+K*w =0,
	\end{equation}
where $r(\phi)=\sum_{n=2}^\infty \frac{F^{(n)}(0)}{n!}\phi^n$.
Since $w$ is odd with respect to $\lambda$ and $K$ is symmetric, we find that
\[
K*w(\lambda)=0
\]
and thus $\beta w_{xx}(\lambda)=0$. If $\beta \neq 0$, then $w$ is smooth (cf. Proposition \ref{prop:beta}), and by Taylor expansion, we have in a neighborhood of $x= \lambda$ that
\[
w(x)=\frac{1}{6}w_{xxx}(\lambda) (x-\lambda)^3 + \frac{1}{24}w_{xxxx}(\xi)(x-\lambda)^4\qquad \mbox{for some}\quad \xi \in [\lambda, x],
\]
where we used that $w(\lambda)=w^{\prime\prime}(\lambda)=0$ and our assumption $w^\prime(\lambda)=0$. In view that $w\geq 0$ on $[\lambda, \lambda+\pi]$, we find from the above expansion that $w_{xxx}(\lambda)\geq 0$. 
Since  $w$ is continuously differentiable in a neighborhood of $\lambda$  , we can take the derivative of \eqref{eq:equation for boundary point lemma 0} in a small neighborhood of $\lambda$ to get
\begin{equation}\label{equation for boundary point lemma1}
	-c w^{\prime}+r^{\prime}(\phi)\phi^\prime-r^{\prime}(\bar\phi)\tilde \phi^\prime +\beta w^{\prime\prime\prime}+[K*w]^{\prime}=0.
\end{equation}
Evaluating the above equation at $x=\lambda$ and using that $r^\prime(\phi)(\lambda)=r^\prime(\tilde \phi)(\lambda)$ together with our assumption  $w^\prime(\lambda)=0$ gives 
\begin{equation}\label{eq:eqc}
[K*w]^\prime(\lambda)=-\beta w^{\prime\prime\prime}(\lambda)\leq 0.
\end{equation}
However,
\begin{align*}
	[K*w]^\prime(\lambda)&= \int_0^\pi K^\prime(y)w(\lambda-y)\,dy + \int_{-\pi}^0 K^\prime(y)w(\lambda-y)\,dy
	=2\int_{-\pi}^0 K^\prime(y)w(\lambda-y)\,dy,
\end{align*}
where we used that $K^\prime$ is antisymmetric and $w$ is odd with respect to $\lambda$. Since $K^\prime>0$ on $(-\pi,0)$ and $w\geq 0$ on $[\lambda,\lambda +\pi]$ but is not identically zero, we obtain that $[K*w]^\prime(\lambda)>0$. This then leads to a contradiction to \eqref{eq:eqc}. Thus, we conclude that either $w$ is identically zero on $\T$ or $w^\prime(\lambda)>0$, which proves the statement.
\end{proof}

{\color{black}As shown in Proposition \ref{prop:beta0}, if $\beta=0$  is not \emph{guaranteed}, then $\phi$ may exhibit a singularity in the form of a peak ($a=2$) or a cusp ($a>2$) so that $\phi$ is merely Lipschitz- or H\"older continuous, respectively.} These cases do not fit into the form of the Hopf boundary point lemma because the $C^{1}$-regularity is not satisfied. However, we can still establish a boundary point lemma for wave profiles with peak or cusp type singularity.

\begin{lem}[Boundary point lemma for highest waves]\label{lem:boundaryH}
	Let $\beta =0$ and  $\phi, \tilde \phi$ be two $2\pi $-periodic continuous solutions of \eqref{eq:GOs} with $B_\phi = B_{\tilde \phi}$. Suppose that $\phi \geq \tilde \phi$ on $[\lambda, \lambda + \pi ]$, $\phi-\tilde \phi$ is odd with respect to $\lambda$. If $\phi(\lambda)=\tilde \phi (\lambda)= \bar \phi$, where $\bar \phi$ is as in \eqref{eq:bar}, then
	\[
	\phi = \tilde \phi \qquad \mbox{on}\quad \T.
	\]
\end{lem}

\begin{proof}
	If $\beta=0$ and $\phi, \tilde \phi$ are two continuous solutions of \eqref{eq:GOs} with $B_\phi = B_{\tilde \phi}$, then
	\[
	F(\phi)- F(\tilde \phi)=-K*w, 
	\]
	where $w=\phi-\tilde \phi$. By Taylor expansion of $F$ around $\bar \phi$, the above equation can be written as
	\[
	\sum_{n\geq a}\frac{F^{(n)}(\bar \phi)}{n!} \left[(\phi-\bar \phi)^n - (\tilde \phi -\bar \phi)^n\right] = -K*w.
	\]
	Here again $a= \min \{ n\geq 2 \mid F^{(n)}(\bar \phi)\neq 0\}$. Since $K*w$ belongs to $C^2(\T)$, we take a Taylor expansion on the right-hand side above around $x=\lambda$ and use the fact $K*w(\lambda)=0$  to obtain
		\[
	\sum_{n\geq a}\frac{F^{(n)}(\bar \phi)}{n!} \left[(\phi-\bar \phi)^n - (\tilde \phi -\bar \phi)^n\right] (x)= -	[K*w]^\prime(\lambda)(x-\lambda) - \frac{1}{2} 	[K*w]^{\prime \prime}(\xi)(x-\lambda)^2\qquad \mbox{for some}\quad \xi \in [\lambda, x].
	\]
 If $w$ is not identically zero, then as in the proof of Lemma \ref{boundary lemma}, we find that 
	\[
	[K*w]^\prime(\lambda)>0
	\]
and
\begin{equation}\label{eq:limit}
\lim_{x\to \lambda}	\sum_{n\geq a}\frac{F^{(n)}(\bar \phi)}{n!} \left[\frac{(\phi-\bar \phi)^n}{x-\lambda} - \frac{(\tilde \phi -\bar \phi)^n}{x-\lambda}\right] (x) = - [K*w]^\prime(\lambda)>0.
\end{equation}
Recall from Proposition \ref{prop:beta0} that $|\phi(x)-\bar \phi|^a \eqsim |x-\lambda|^2$ and the same is true for $\tilde \phi$. 
So,  the limit on the left-hand side of \eqref{eq:limit} equals zero, which leads to a contradiction.
\end{proof}

With the above lemmata, we are now ready to prove our main result on the symmetry of periodic solutions of \eqref{eq:GOs}. {\color{black} We distinguish the cases $\beta=0$ and $\beta>0$. In the case $\beta=0$, the traveling wave solution might not be smooth but exhibit a singularity at the crest or trough. In this case, our symmetry result requires that the solution has a \emph{unique} global maximum and minimum per period, but there are no further monotonicity assumptions on the wave profile. Then, we prove in Theorem \ref{thm:symmetry_beta_0} that the traveling wave is symmetric with a single crest per period. In the case $\beta>0$, the traveling wave solution is smooth, but we have to impose a priori a  monotonicity assumption on the wave profile in order to prove the symmetry result in Theorem \ref{thm:symmetry_beta_1}.}
\begin{thm}[Symmetry of periodic solutions for $\beta=0$] \label{thm:symmetry_beta_0}
	Let $\phi$ be a nontrivial, continuous $2\pi$-periodic solution of \eqref{eq:GOs} with $\beta=0$ and $F^\prime(\phi)\leq 0$. If $\phi$ has a unique global maximum and minimum in $\T$, then $\phi$ is symmetric. Moreover,
	\[
	\phi^\prime(x)>0\qquad \mbox{for all}\qquad x\in (-\pi,0),
	\]
	after a possible translation.
\end{thm}
\begin{proof}
The condition $F^\prime(\phi)\leq 0$  implies that the solution $\phi$ satisfies
	\[
	\max \phi \leq \bar \phi_+ \qquad \mbox{and}\qquad \min \phi \geq \bar \phi_-,
	\]
	whenever $\bar \phi_+$ or $\bar \phi_-$ exists. Recall that $\bar \phi_+>0$ ($\bar \phi_-<0$) is the smallest positive (largest negative) solution of $F^\prime(\phi)=0$ if it exists. 
From Remark \ref{rem:amplitude} we know that $F^\prime(\phi)\leq 0$ is a priori satisfied if $a=\min \{n\geq 2 \mid F^{(n)}(\bar \phi)\neq 0\}$ is even. 
Assume without loss of generality that $\phi$ is a nontrivial solution and that the unique global minimum in $\T$ is located at $x=-\pi$, that is $\min \phi =\phi(-\pi)=\phi(\pi)$. Moreover, suppose that $\phi$ attains its unique global maximum at $x_M\in \T$.\\
	\emph{Monotonicity}:
		In order to show the one-sided monotonicity of the wave profile, we are going to use the local formulation of \eqref{eq:GOs} with $\beta=0$, that is, $\phi$ is a continuous function, solving
	\begin{equation}\label{eq:lf}
		(F(\phi))_{xx}=\phi
	\end{equation}
pointwise almost everywhere.
 First, we use a modification of the so-called \emph{sliding method} (cf. \cite{BN91}) to prove the monotonicity of $\phi$ in the region $[-\pi,x_M]$ between the unique global trough and the unique global maximal crest. It is clear that $\phi (-\pi)\leq \phi(x)<\phi(x_M)$ for all $x\in (-\pi, x_M)$. If $\phi$ is a solution of \eqref{eq:lf}, then by translation invariance also $\phi_h:=\phi(\cdot+h)$ is a solution for any $h\in \R$. Note that $\phi_{x_M+\pi}$ attains its global maximum at $x=-\pi$. For $h\in [0,x_M+\pi]$, we set 
 \[\Omega_h:=[-\pi,x_M-h]
 \] 
 and
 \[
 h_*:=\inf\{h\in [0,x_M+\pi]\mid \phi_{h}-\phi \geq 0\quad \mbox{on}\quad \Omega_h\}.
 \]
Clearly, such an infimum exists, since $\Omega_{x_M}=\{-\pi\}$ and $\phi_{x_M}-\phi>0$ at $x=-\pi$. We aim to show that $h_*=0$, that is, $\phi$ is monotone on $\Omega_0$.

\begin{figure}[h!]
\begin{tikzpicture}[scale=1,yscale=1.3]
	\small
	\draw (-2,0)--(8,0);
	\draw[-,black](0,0.1)--(0,-0.1) node[below] {$-\pi$};
	\draw[-,black](6,0.1)--(6,-0.1) node[below] {$\pi$};
	\draw[luh-dark-blue, very thick] plot [smooth] coordinates {(-1,1.3) (0,0.5) (0.8,1) (1.3,0.8) (2.1,1.4) (2.5, 1.8)};
		\draw[luh-dark-blue, very thick] plot [smooth] coordinates { (2.5, 1.8) (3.2, 1.6) (5.1,1.4)  (6,0.5) (6.8,0.8)} node[right]{$\phi$};
	\begin{scope}[xshift=-1.2cm]
		\draw[gray,  thick] plot [smooth] coordinates {(-1,1.3) (0,0.5) (0.8,1) (1.3,0.8) (2.1,1.4) (2.5, 1.8)} node[left]{$\phi_h$};
	\draw[gray,  thick] plot [smooth] coordinates { (2.5, 1.8) (3.2, 1.6) (5.1,1.4)  (6,0.5) (6.8,0.8)};
	\end{scope}
\draw[orange, very thick] (0,0)--(1.3,0);
\node at (0.7,0.2) {\textcolor{orange}{$\Omega_{h}$}};
	\draw[-] (2.5,0.1)--(2.5,-0.1)node[below] {$x_M$};
\end{tikzpicture}
\caption*{Illustration of the sliding method.}
\end{figure}
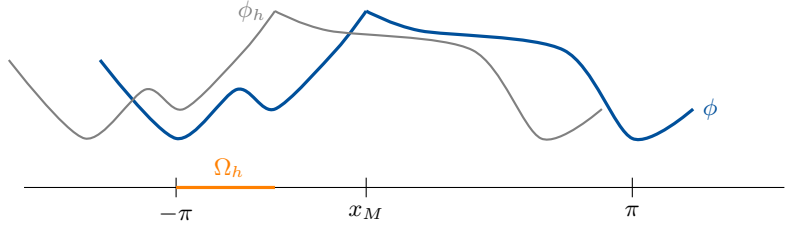

To proceed, we assume for a contradiction that $h_*>0$ and there exists $\bar x \in \Omega_{h_*}$ such that $\phi_{h_*}(\bar x)=\phi(\bar x)$ and $\phi_{h_*}^\prime(\bar x)=\phi^\prime(\bar x)$. Since $\phi$ is assumed to have a unique global maximum and minimum in $\T$,  we find that $(\phi_{h_*}-\phi) \big|_{x\in \partial \Omega_{h_*}}>0$ so that that $\bar x$ is an inner point of $\Omega_{h_*}$, that is $\bar x \in \Omega_{h_*}^\circ:=(-\pi, x_M-h_*)$. {\color{black}Notice that $\phi$ and $\phi_h$ are smooth on $\Omega_{h_*}^\circ$ and satisfy \eqref{eq:lf} pointwise.}
Setting $v_{h_*}:=\phi_{h_*}-\phi$ and using assumptions (F) we find that
{
 \begin{equation}\label{eq:lf1}
 (F(\phi)-F(\phi_{h_*}))_{xx}=-v_{h_*},
 \end{equation}
 }
 Denote by $f:=F(\phi)-F(\phi_{h_*})$. Then $f$ is twice continuously differentiable with
 {
 \[
 f^{\prime\prime} \leq  0 \qquad \mbox{on}\quad \Omega_{h_*}
 \]
 } 
and  we use the mean value theorem to get
 \[
 f=F^\prime(\phi^\star)(\phi-\phi_{h_*})\qquad \mbox{where} \quad  \phi^\star \in [\phi, \phi_{h_*}].
 \]
Since $\phi( x)\neq \bar \phi$ and $\phi_{h_*}( x)\neq \bar \phi$ on $\Omega_{h_*}^\circ$, where $\bar \phi \in \{\bar \phi_+, \bar \phi_-\}$, it is clear from the assumptions (F) that 
 \[
 F^\prime(\phi^\star)(x)<0\qquad \mbox{for all}\quad x \in\Omega_{h_*}^\circ.
 \]
 In particular, this indicates { $f \geq 0$} on $\Omega_{h_*}$ with $f(\bar x)=0$ and $\bar x \in \Omega_{h_*}^\circ$. Since $f$ is {concave} on $\Omega_{h_*}^\circ$, we conclude by  maximum principle  that  $f=0$ on $\Omega_{h_*}$. However,  this is only possible if $h_*=0$, which contradicts the assumption that $h_*>0$. Hence, we conclude that $\phi$ is monotone (increasing) on $[-\pi,x_M]$.
In view that \eqref{eq:lf} is also invariant under symmetric reflection, we can use the same argument on the reflected function $\tilde \phi:=\phi(-\cdot)$ to prove that $\phi$ is also monotone (decreasing) on $[x_M,\pi]$. \\
 \emph{Symmetry:} Having established that $\phi$ has a single crest per period, we now use the method of moving planes to show that $\phi$ is symmetric and $\phi^\prime(x)>0$ on $(-\pi,0)$. Set
 \[
 w_\lambda:=\phi(x)-\phi(2\lambda-x).
 \]
 Note that if $\phi$ is a solution of \eqref{eq:GOs}, then the reflection about $\lambda$ given by $\phi(2\lambda-\cdot)$ is a solution of $\eqref{eq:GOs}$ due to the symmetry of $K$. Since $\phi$ is monotone on $(-\pi,x_M)$ and $2\pi$-periodic, there exists $\lambda_0\in (-\pi,x_M)$ such that $w_{\lambda_0}>0$ on $(\lambda_0,\lambda_0+\pi)$. 
 	Set 
 \[\lambda_*:= \sup \{\lambda \in [\lambda_0,0 ] \mid w_\lambda(x) > 0 \; \mbox{ for all } \; x\in (\lambda,\lambda+\pi )\}.\]
 We have that $w_{\lambda_*}$ satisfies
 \begin{itemize}
 	\item[i)] $w_{\lambda_*}({\lambda_*})=0$; 
 	\item[ii)] $w_{\lambda_*}$ is odd with respect to $\lambda_*$, that is $w_{\lambda_*}(\cdot)=-w_{\lambda_*}(2\lambda_*-\cdot)$;
 	\item[iii)] $w_{\lambda_*}\geq 0$ in $[\lambda_*, \lambda_*+\pi ]$ and $w_{\lambda_*}\leq 0$ in $[\lambda_*+\pi , \lambda_*+2\pi ]$.
 \end{itemize}
If $w_{\lambda_*}$ is identically zero, then $\phi$ is symmetric. So, we assume now for a contradiction that $w_{\lambda_*}$ is not identically zero. We start with $\lambda=\lambda_0$, reflect the  wave profile about the axis $x=\lambda$, and then push forward the axis $x=\lambda$ as long as $w_\lambda>0$ on $(\lambda, \lambda+\pi)$. This procedure will stop at $\lambda_*$. 
 	Clearly, the process stops at or before the crest in $x=x_M$ at $\lambda=\lambda_*$. In fact, one of three occasions will occur: Either there exists $\bar x \in (\lambda_*,\lambda_*+\pi )$ such that $w_{\lambda_*}(\bar x)=0$; or we reach a crest at $x=\lambda_*$; or we reach a  trough at $\lambda_*+\pi $.
 The first case can be excluded by the touching lemma (cf. Lemma \ref{touching lemma}) unless $w_{\lambda_*}$ is identically zero. 
 For the latter two cases, if $\phi$ is continuously differentiable in a neighborhood of $\lambda_*$, we can apply  Lemma \ref{boundary lemma}  to obtain the contradiction $\phi^\prime(\lambda)>0$, while $\phi$ reaching a local maximum or minimum at $\lambda=\lambda_*$ or $\lambda=\lambda_* + \pi$, respectively. If $\phi$ is a highest wave with $\phi(\lambda_*)= \bar \phi_+$ or $\phi(\lambda_* + \pi)=\bar \phi_-$, then Lemma \ref{lem:boundaryH} implies that $w_{\lambda_*}$ is identically zero. 
Therefore, for all these three cases, we conclude that $\phi$ is symmetric and has exactly one crest per period. 

Since we fix the trough at $x=-\pi$, we deduce from the symmetry of $\phi$ that the crest is located at $x=0$. Moreover, by the assumption that $\phi$ has a unique global maximum and minimum and $F^\prime(\phi)\leq 0$, we actually have that $F^\prime(\phi)<0$ on $(-\pi, 0)$. Then, the fact that $\phi$ has exactly one crest per period and Proposition \ref{prop:beta0} imply that  $\phi^\prime(x)\geq 0$ for all $x\in (-\pi,0)$.
 We are left to show that the strict inequality prevails for any $x\in (-\pi,0)$.  Taking derivative on equation \eqref{eq:rGOt}, we obtain for $x\in (-\pi,0)$ that
 \[
 F^\prime(\phi)\phi^{\prime}(x)=-[K*\phi]^\prime (x)= -\int_{-\pi}^0 \left[ K(x-y)-K(x+y) \right]\phi^\prime(y)\,dy,
 \]
 where we used the symmetry of $K$ and $\phi$. In the same fashion as in the proof of the touching lemma (cf. Lemma \ref{touching lemma}), one can show that the right-hand side of the above equation is strictly negative unless $\phi$ is a trivial solution. In view of $F^\prime(\phi)<0$ on $(-\pi,0)$, we conclude the proof.
\end{proof}

{\color{black} In the case when $\beta>0$, we can not repeat the same argument by sliding method as in the proof of the previous theorem, since \eqref{eq:lf} would become
\[
(F(\phi)+\beta \phi_{xx})_{xx}=\phi
\]
and the corresponding equation for the difference of $\phi$ and its shift $\phi_{h_*}$ in \eqref{eq:lf1} reads
\[
(F(\phi)-F(\phi_{h_*})-\beta (v_{h_*})_{xx})_{xx}=-v_{h_*}.
\]
We could still conclude that $f:=F(\phi)-F(\phi_{h_*})+\beta (v_{h_*})_{xx}=-F^\prime(\phi^\star)v_{h_*}-\beta (v_{h_*})_{xx}$ is concave on $\Omega_{h_*}$, but we don't have enough information about the sign of $f$ in order to force a contradiction. Therefore, we impose in the following theorem on the symmetry of traveling waves for $\beta>0$, that the wave profile has a priori a single crest per period.
}

\begin{thm}[Symmetry of periodic traveling waves for $\beta >0$]\label{thm:symmetry_beta_1}
	Let $\phi$ be a nontrivial, $2\pi$-periodic solution of \eqref{eq:GOs}. If $\phi$ has a single crest per period, then $\phi$ is symmetric. 
\end{thm}

\begin{proof}
	Recall that   $\phi$ is smooth for  $\beta>0$ (cf. Proposition \ref{prop:beta})). Since we assume a priori that $\phi$ has a single crest per period, the proof follows the lines of the second part (symmetry) in the proof of Theorem \ref{thm:symmetry_beta_0} by using Lemma \ref{touching lemma} and the boundary point lemma for continuously differentiable functions (Lemma \ref{boundary lemma})).

\end{proof}

\begin{remark}\emph{
The (quite strong) assumption in Theorem \ref{thm:symmetry_beta_1} that the wave profile has a single crest per period can be relaxed to the assumption	that there exists $\lambda_0\in [0,2\pi)$ such that
\[
{\color{black}\phi(x)>\phi(2\lambda_0-x)}\qquad \mbox{for all}\qquad x\in (\lambda_0,\lambda_0+\pi).
\]
 This condition is referred to as a \emph{reflection criterion} in \cite{2021Symmetry}, where it is shown that the method of moving planes can be applied to prove the symmetry of periodic traveling waves for a class of equations with very weak dispersion. With small adaptations, the same argument can be applied in the case at hand.
}
\end{remark}

\bigskip

\section{Application to the Gardner--Ostrovsky equation}\label{sect:application to GO}

The Gardner--Ostrovsky equation is given by
\begin{equation}\label{GOE}
(u_t + \sigma uu_x + \alpha u^2u_x+\beta u_{xxx})_x = u
\end{equation}
and $u(t,x)=\phi(x-ct)$ is a traveling wave solution with wave speed $c>0$ if $\phi$ satisfies the 
 steady Gardner--Ostrovsky equation 
\begin{equation}\label{GOES}
	F(\phi)+\beta \phi_{xx}=-D^{-2}\phi+B_\phi,
\end{equation}
where 
\[
F(\phi)=-c \phi + \frac{\sigma}{2}\phi^2 + \frac{\alpha^3}{3}\phi^3.
\]
Here  $\alpha, \sigma \in \R$ are free parameters. For $\alpha=0$, we recover the Ostrovsky equation, and for $\sigma=0$, the modified Ostrovsky equation (also called the short pulse equation). 
The following corollaries summarize our previously established results {\color{black}about the existence, regularity and amplitude of periodic traveling solutions} when applied to the Gardner--Ostrovsky equation. 

\begin{cor}
 For $\beta \in [0,1)$ there exists a local {\color{black}smooth} bifurcation curve of nontrivial $2\pi$-periodic real-valued solutions of \eqref{GOES} emanating from the bifurcation point $(c_*,0)\in \R \times C^a_{0,\textrm{even}}(\T)$, where $c_*=1-\beta$ and $a\in (1,2)$. In addition, any continuous, $2\pi$-periodic solution $\phi$ satisfies 

	\begin{itemize}
	\item If $\beta>0$, then  $\phi$ is smooth, and $\phi$ is symmetric if it has a single crest per period.
	\item If $\beta=0$, then $\phi$ is smooth on any open set where $\sigma \phi + \alpha \phi^2<c$. In particular, $\phi$ is smooth if $\alpha<-\frac{\sigma^2}{4c}$. However, 
	$\phi$ does not belong to $C^1(\T)$ if $\phi$ has an isolated maximum or minimum at $\bar x$ with $\phi(\bar x)=\bar \phi$ such that \[\sigma \bar \phi + \alpha \bar \phi^2=c.\] More precisely, $\phi$ does not belong to $C^1(\T)$ in the following cases:\\[5pt]
	\begin{center}
	\begin{tabular}{ |c | c| c| c| c| c }
				\hline
		& $\alpha =- \frac{\sigma^2}{4c}$ & $- \frac{\sigma^2}{4c}<\alpha<0$ & $\alpha=0$& $\alpha>0$ \\ 
		\hline
	$\sigma>0$&	$\max \phi = \phi_*^-$ & $\max \phi = \phi_*^-$ & $\max \phi=\phi_{*,0}$& 	$\max \phi = \phi_*^+$ or 	$\min \phi = \phi_*^-$  \\  	\hline
	$\sigma=0$&	$\times$ & 	$\times$ & 	$\times$ &	$\max \phi = \phi_*^+$ or 	$\min \phi = \phi_*^-$ \\		\hline
	$\sigma<0$  &$\min \phi = \phi_*^+$ &$\min \phi = \phi_*^+$ & $\min \phi=\phi_{*,0}$&	$\max \phi = \phi_*^+$ or 	$\min \phi = \phi_*^-$ \\		\hline
	\end{tabular}
\end{center}

\vspace{10pt}

\begin{itemize}
	\item[]
Here, $\phi_*^\pm$ are given by
\[
	\phi_*^\pm =-\frac{\sigma}{2\alpha}\pm \sqrt{\frac{\sigma^{2}}{4\alpha^2}+\frac{c}{\alpha}}, \qquad \qquad \phi_{*,0}=\frac{c}{\sigma}.
\]
If $\alpha>-\frac{\sigma^2}{4c}$ and $\phi$ satisfies the condition in the table above at $x=0$\footnote{Without loss of generality, one can shift the wave profile such that the condition in the table is satisfied by $\phi$ at $x=0$ due to the translation invariance of  equation \eqref{GOES}.}, then $\phi$ is precisely Lipschitz continuous at $x=0$ with
\[
|\phi^\prime(0_{\pm})|^2= \frac{1}{2|\sigma + 2\alpha \phi(0)|}|\phi(0)|.
\]
If $\alpha=-\frac{\sigma^2}{4c}$ and $\phi$ satisfies the condition in the table above at $x=0$, then $\phi$ is precisely $\frac{2}{3}$-H\"older continuous at $x=0$ with
\[
|\phi (x)-\phi(0)|\eqsim \left( \frac{3}{2|\alpha|}\right)^{\frac{1}{3}}|x|^{\frac{2}{3}}\qquad \mbox{for}\quad |x|\ll 1.
\]

 \end{itemize}
\item  The amplitude of $\phi$ is characterized by the following:
 \begin{itemize}
	\item[]
 If $\alpha >0$, then the amplitude of $\phi$ is bounded by
\[
\max \phi - \min \phi \leq 2 \sqrt{\frac{\sigma^2}{4\alpha^2}+\frac{c}{\alpha}}.
\]
\item []
If $-\frac{\sigma^2}{4c}<\alpha\leq 0$, then $\phi$ is bounded from above by $\phi_*^->0$ ($\phi_{*,0}$ for $\alpha=0$) if $\sigma>0$ and $\phi$ is bounded from below  by $\phi_*^+<0$ ($\phi_{*,0}$ for $\alpha=0$) if $\sigma<0$. 
\end{itemize}
\end{itemize}

\end{cor}

The  highest waves are explicitly given in terms of polynomials of degree $2$ for the reduced Ostrovsky family and of degree $1$ for the modified reduced Ostrovsky family (see \cite{MR3685174} and \cite{BD}) as follows:
\begin{itemize}
	\item 
The highest wave for the reduced Ostrovsky family with $\alpha=0$ and $\sigma\in \R\setminus\{0\}$ is given by
\[
\phi_\sigma(x)= \frac{3(|x|-\pi^2)^2-\pi^2}{18\sigma} \qquad \mbox{with wave speed}\qquad c=\frac{\pi^2}{9}.
\]
\item 
The highest wave for the modified reduced Ostrovsky family with $\sigma=0$ and $\alpha>0$ is given by
\[
\phi_\alpha(x)= \frac{1}{\sqrt{2\alpha}}\left(\frac{\pi^2}{2}-|x|\right) \qquad \mbox{with wave speed}\qquad c=\frac{\pi^2}{8}.
\]
Observe that the parameters $\sigma$ and $\alpha$ above do not influence the speed of the highest wave but only on its amplitude (and angle enclosed by the peaked singularity). 

\begin{center}
	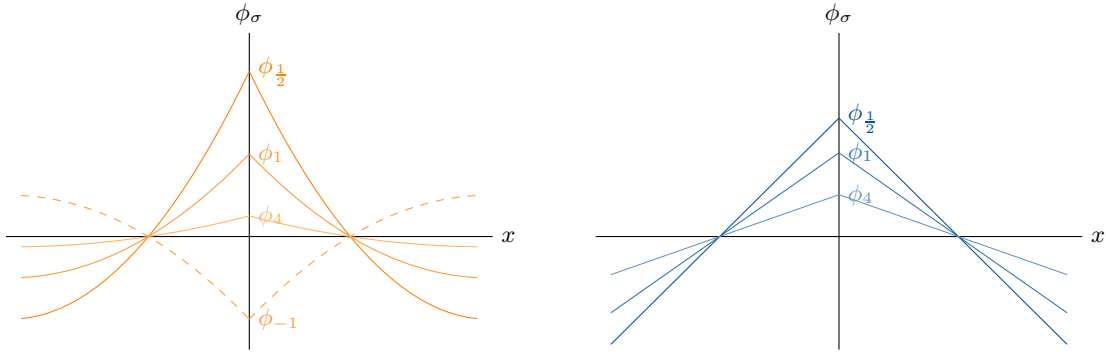
\begin{figure}[h!]
\begin{tikzpicture}
	\small
	\draw[-] (-3.2,0)--(3.2,0) node[right] {$x$} ;
	\draw[-] (0,-1.5)--(0,2.7) node[above] {$\phi_{\sigma}$};
\draw[color=orange!80, domain=0:3] plot(\x,{1/18 *(3*(\x-3.14)*(\x-3.14)-3.14*3.14)}) ;
\draw[color=orange!80, domain=-3:0] plot(\x,{1/18 *(3*(-\x-3.14)*(-\x-3.14)-3.14*3.14)})node[right] {$\phi_1$} ;

\draw[color=orange, domain=0:3] plot(\x,{1/(18*0.5) *(3*(\x-3.14)*(\x-3.14)-3.14*3.14)}) ;
\draw[color=orange, domain=-3:0] plot(\x,{1/(18*0.5) *(3*(-\x-3.14)*(-\x-3.14)-3.14*3.14)}) node[right] {$\phi_{\frac{1}{2}}$} ;

\draw[color=orange!60, domain=0:3] plot(\x,{1/(18*4) *(3*(\x-3.14)*(\x-3.14)-3.14*3.14)}) ;
\draw[color=orange!60, domain=-3:0] plot(\x,{1/(18*4) *(3*(-\x-3.14)*(-\x-3.14)-3.14*3.14)})node[right] {$\phi_4$} ;

\draw[color=orange!80, dashed, domain=0:3] plot(\x,{1/(18*(-1)) *(3*(\x-3.14)*(\x-3.14)-3.14*3.14)}) ;
\draw[color=orange!80, dashed, domain=-3:0] plot(\x,{1/(18*(-1)) *(3*(-\x-3.14)*(-\x-3.14)-3.14*3.14)}) node[right] {$\phi_{-1}$};
\end{tikzpicture}
\hspace{2em}
	\begin{tikzpicture}
		\small
		\draw[-] (-3.2,0)--(3.2,0) node[right] {$x$};
		\draw[-] (0,-1.5)--(0,2.7)  node[above] {$\phi_{\sigma}$};
		\draw[color=luh-dark-blue!80, domain=0:3] plot(\x,{(1/sqrt(2*1)) *(3.14*0.5-\x)});
		\draw[color=luh-dark-blue!80, domain=-3:0] plot(\x,{(1/sqrt(2*1)) *(3.14*0.5+\x)})node[right] {$\phi_1$}  ;
		
		\draw[color=luh-dark-blue, domain=0:3] plot(\x,{(1/sqrt(2*0.5)) *(3.14*0.5-\x)}) ;
	\draw[color=luh-dark-blue, domain=-3:0] plot(\x,{(1/sqrt(2*0.5)) *(3.14*0.5+\x)}) node[right] {$\phi_{\frac{1}{2}}$} ;
		
		\draw[color=luh-dark-blue!60, domain=0:3] plot(\x,{(1/sqrt(2*4)) *(3.14*0.5-\x)}) ;
	\draw[color=luh-dark-blue!60, domain=-3:0] plot(\x,{(1/sqrt(2*4)) *(3.14*0.5+\x)})node[right] {$\phi_4$}  ;
	\end{tikzpicture}
\caption{Plots of highest waves for the reduced Ostrovsky equation ($\alpha=0$) on the left and for the modified reduced Ostrovsky equation ($\sigma=0$) on the right.} 
\end{figure}
\end{center}
\end{itemize}

\subsection*{Acknowledgment}
The author L.P. gratefully
acknowledges financial support from the National Natural Science Foundation for Young Scientists of China
(Grant No. 12001553), the Fundamental Research Funds for the Central Universities (Grant No. 20lgpy151),  the Science and Technology Program of Guangzhou (Grant
No. 202102080474) and the {\color{black} Guangdong Basic,   and Applied Basic Research Foundation (Grant No. 2023A1515010599)}.

\subsection*{Statements and Declarations}
Competing Interest: this work does not have any conflicts of interest.

%
%

	\end{document}